\documentclass[10pt,final,twocolumn]{IEEEtran}
\long\def\comment#1{}

\usepackage{textcomp}
\usepackage{cite}
\usepackage{graphicx}
\usepackage{subfigure}
\usepackage{psfrag}
\usepackage{url}
\usepackage{stfloats}
\usepackage{amsmath}
\usepackage{amssymb,amsfonts}
\usepackage{amsthm}
\usepackage{float}
\usepackage[usenames]{color}
\usepackage{algorithm,algorithmic}

\newtheorem{assumption}{Assumption}
\newtheorem{remark}{Remark}

\newtheorem{theorem}{Theorem}
\newtheorem{example}{Example}
\newtheorem{proposition}{Proposition}

\begin{document}

\setlength{\arraycolsep}{0.3em}

\title{On Faster Convergence of Scaled Sign Gradient Descent
\thanks{}}

\author{Xiuxian Li, Kuo-Yi Lin, Li Li, Yiguang Hong, and Jie Chen
\thanks{This work was supported by the Shanghai Municipal Commission of Science and Technology No. 19511132100, 19511132101, the Shanghai Municipal
Science and Technology Major Project, No. 2021SHZDZX0100, and the National Natural Science Foundation of China under Grant 62003243.}
\thanks{The authors are with Department of Control Science and Engineering, College of Electronics and Information Engineering, Shanghai Research Institute for Intelligent Autonomous Systems, and Shanghai Institute of Intelligent Science and Technology, Tongji University, Shanghai, China (e-mail: xli@tongji.edu.cn, 19603@tongji.edu.cn, lili@tongji.edu.cn, yghong@tongji.edu.cn, chenjie206@tongji.edu.cn).}
}

\maketitle

\setcounter{equation}{0}
\setcounter{figure}{0}
\setcounter{table}{0}

\begin{abstract}
Communication has been seen as a significant bottleneck in industrial applications over large-scale networks. To alleviate the communication burden, sign-based optimization algorithms have gained popularity recently in both industrial and academic communities, which is shown to be closely related to adaptive gradient methods, such as Adam. Along this line, this paper investigates faster convergence for a variant of sign-based gradient descent, called scaled {\scriptsize SIGN}GD, in three cases: 1) the objective function is strongly convex; 2) the objective function is nonconvex but satisfies the Polyak-{\L}ojasiewicz (PL) inequality; 3) the gradient is stochastic, called scaled {\scriptsize SIGN}SGD in this case. For the first two cases, it can be shown that the scaled {\scriptsize SIGN}GD converges at a linear rate. For case 3), the algorithm is shown to converge linearly to a neighborhood of the optimal value when a constant learning rate is employed, and the algorithm converges at a rate of $O(1/k)$ when using a diminishing learning rate, where $k$ is the iteration number. The results are also extended to the distributed setting by majority vote in a parameter-server framework. Finally, numerical experiments on logistic regression are performed to corroborate the theoretical findings.
\end{abstract}


\section{Introduction}\label{s1}

This paper studies an unconstrained optimization problem
\begin{align}
\min_{x\in\mathbb{R}^d} f(x),             \label{1}
\end{align}
where the objective $f:\mathbb{R}^d\to\mathbb{R}$ is a proper differentiable function, and may be nonconvex, which has numerous applications in industry, such as electric vehicles \cite{shen2014optimization,li2018distributedon}, smart grid \cite{su2011survey}, internet of things (IoT) \cite{messaoud2019online}, and so on. To solve this problem, a quintessential algorithm is the gradient descent (GD) method \cite{ruder2016overview,meng2020aug,li2021distributed}, which requires to access true gradients. However, it is usually expensive or difficult to compute the true gradients in reality, and thereby a typical stochastic gradient descent (SGD) algorithm has become prevalent in deep neural networks \cite{bottou2012stochastic,bonnabel2013stochastic}, which depends upon a lower computing cost for stochastic gradients.

As for large-scale neural networks, the training efficiency can be substantially improved in general by introducing multiple workers in a parameter-server framework, where a group of workers can train their own mini-batch datasets in parallel. Nonetheless, the communication between workers and the parameter server has been a non-negligible handicap for its wide practical application. As such, as one of gradient compression techniques, sign-based methods have been popular in recent decades, not only because they can reduce the communication cost to one bit for each gradient coordinate, but because they have good performance and close relationship with adaptive gradient methods \cite{riedmiller1993direct,bernstein2018signsgd,balles2020geometry}. As a matter of fact, it has been demonstrated in \cite{balles2018dissecting,bernstein2018signsgd} that {\scriptsize SIGN}SGD with momentum often has pretty similar performance to Adam on deep learning missions in practice. Notice that a wide range of gradient compression approaches exist for reducing the communication cost in the literature, e.g., \cite{hamer2020fedboost,xie2020cser}, whose elaboration is beyond the scope of this paper. Particularly, sign-based methods considered in this paper can be regarded as a special gradient compression scheme which need to transmit only one bit per gradient component \cite{karimireddy2019error}.

Along this line, the sign gradient descent ({\scriptsize SIGN}GD) algorithm and its stochastic counterpart ({\scriptsize SIGN}SGD) have been extensively studied in recent years \cite{seide20141,bernstein2018signsgd,bernstein2019signsgd,balles2020geometry,safaryan2019stochastic}, which are, respectively, of the form
\begin{align}
x_{k+1}&=x_k-\alpha_k \operatorname{sign} (\nabla f(x_k)),        \label{2}\\
x_{k+1}&=x_k-\alpha_k \operatorname{sign} (h_k),                  \label{3}
\end{align}
where $h_k\in\mathbb{R}^d$ is a stochastic gradient of $f$ at $x_k$, $\alpha_k>0$ is the learning rate, and the signum function $\operatorname{sign}$ is operated componentwise. For instance, it was demonstrated in \cite{bernstein2018signsgd} that {\scriptsize SIGN}SGD enjoys a SGD-level convergence rate for nonconvex but smooth objective functions under a separable smoothness assumption, which, in combination with majority vote in distributed setup, was further shown to be efficient in terms of communication and fault toleration in \cite{bernstein2019signsgd}. Recently, the authors in \cite{balles2020geometry} found that the $\ell_\infty$-smoothness is a weaker and natural assumption than the separable smoothness and established two conditions under which the sign-based methods are preferable over GD.

{\bf Contributions.}
To our best knowledge, this paper is the first to address faster convergence of sign methods with more details as follows.


First, it is found that {\scriptsize SIGN}GD is not generally convergent even for strongly convex and smooth objectives when using constant learning rates, although it is indeed convergent for vanilla GD. Therefore, scaled versions in Algorithms \ref{al1} and \ref{al2} are investigated. It is proved that Algorithm \ref{al1} converges linearly to the minimal value for two cases: strongly convex objectives and nonconvex objectives yet satisfying the Polyak-{\L}ojasiewicz (PL) inequality. Meanwhile, Algorithm \ref{al2} converges linearly to a neighborhood of the minimal value when using a constant learning rate $\alpha$ with an error being proportional to $\alpha^2$ and the variance of stochastic gradients. When applying a kind of diminishing learning rate, a rate $O(\ln^2(k)/k^2)$ can be ensured for (\ref{5}), which is superior to the widely known rate $O(1/k)$ \cite{stich2018sparsified}.


Second, the obtained results are extended to the distributed setup, where a group of workers compute their own (stochastic) gradients using individual dataset and then transmit the sign gradient and the gradient $\ell_1$-norm to the parameter server who calculates the sign gradient by majority vote along with taking the average of the gradient $\ell_1$-norms and transmits back to all the workers.


{\bf Notations.}
Denote by $[n]:=\{1,2,\ldots,n\}$ for an integer $n>0$. Let $\|\cdot\|$, $\|\cdot\|_1$, $\|\cdot\|_\infty$ and $x^\top$ be the $\ell_2$-norm, $\ell_1$-norm, $\ell_\infty$-norm and the transpose of $x\in\mathbb{R}^n$, respectively. $\mathbf{1}$ and $\mathbf{0}$ stand for column vectors of compatible dimension with all entries being $1$ and $0$, respectively. $\nabla f$ represents the gradient of a function $f$. $\mathbb{E}(\cdot)$ and $\mathbb{P}(\cdot)$ denote the mathematical expectation and probability, respectively.

\section{Counterexamples for {\footnotesize SIGN}GD}\label{s3}

For {\scriptsize SIGN}GD, an interesting result can also be found in the continuous-time setup, which demonstrates obvious advantages of {\footnotesize SIGN}GD compared with GD. Particularly, {\scriptsize SIGN}GD converges linearly, while GD is only sublinearly convergent. More details are postponed to the Appendix as supplemental materials.

Motivated by the fact in the continuous-time setup, it seems promising to consider the discrete-time counterpart of (\ref{ca2}), i.e., {\scriptsize SIGN}GD (\ref{2}) with $\alpha_k=\alpha>0$ being a constant learning rate. However, it is not the case. It is well known that GD is linearly convergent for small enough $\alpha>0$, while several counterexamples are presented below for illustrating that the sign counterpart (\ref{2}) is generally not convergent even for strongly convex and smooth objectives.

\begin{example}\label{e1}
Consider $f(x)=x_1^2+x_2^2$ for $x\in\mathbb{R}^2$, which is strongly convex and smooth with $\nabla f(x)=(2x_1,2x_2)^\top$. By choosing the initial point as $x_0=(\alpha/2,\alpha/2)^\top$, it is easy to verify for (\ref{2}) that for $l=0,1,2,\ldots$,
\begin{align}
x_{2l}=\big(-\frac{\alpha}{2},-\frac{\alpha}{2}\big)^\top,~~~
x_{2l+1}=\big(\frac{\alpha}{2},\frac{\alpha}{2}\big)^\top,             \label{8}
\end{align}
which is obviously not convergent.
\end{example}

Example \ref{e1} shows that the exact convergence cannot be ensured for {\scriptsize SIGN}GD even for strongly convex and smooth objectives. To fix it, one may attempt to consider the sign counterpart of adaptive gradient methods. However, it generally does not work as well. For instance, the AdaGrad-Norm \cite{ward2019adagrad}
\begin{align}
b_{k+1}^2&=b_k^2+\|\nabla f(x_k)\|^2,         \nonumber\\
x_{k+1}&=x_k-\frac{\eta}{b_{k+1}}\nabla f(x_k),~~~\eta>0        \label{9}
\end{align}
is shown to converge linearly without knowing any function parameters beforehand \cite{xie2020linear}, while the linear convergence cannot be ensured in general for its sign counterparts, as illustrated below for its two sign variants.

\begin{example}\label{e3}
Consider the first sign variant as
\begin{align}
b_{k+1}^2&=b_k^2+\|\nabla f(x_k)\|^2,         \nonumber\\
x_{k+1}&=x_k-\frac{\eta}{b_{k+1}}\operatorname{sign}(\nabla f(x_k)),~~~\eta>0        \label{10}
\end{align}
and $f(x)=x^2/2$ (strongly convex and smooth) with $x\in\mathbb{R}$. For simplicity, set $b_0=0$ and $x_0\neq 0$. Then simple manipulations give rise to $b_{k+1}^2=\sum_{l=0}^k x_l^2$.

In what follows, we show that the convergence rate of (\ref{10}) is not linear. To do so, it is easy to see that $x_{k+1}=x_k-\frac{\eta}{b_{k+1}}\operatorname{sign}(x_k)$, which leads to that
\begin{align}
x_{k+1}^2&=x_k^2-\frac{2\eta}{b_{k+1}}|x_k|+\frac{\eta^2}{b_{k+1}^2}           \nonumber\\
&=x_k^2-\frac{\eta}{b_{k+1}^2}(2b_{k+1}|x_k|-\eta).                        \label{12}
\end{align}
By contradiction, if $x_k$ or $f(x_k)$ is linearly convergent, then one has that $\sum_{k=0}^\infty x_k^2\leq B$ for some constant $B>0$, which, together with (\ref{12}) and $b_{k+1}^2=\sum_{l=0}^k x_l^2$, gives
\begin{align}
x_{k+1}^2\geq x_k^2-\frac{\eta}{b_{k+1}^2}(2\sqrt{B}|x_k|-\eta).          \label{13}
\end{align}
After $|x_k|$ decreases to where $|x_k|\leq \frac{\eta}{2\sqrt{B}}$, invoking (\ref{13}) leads to $x_{k+1}^2\geq x_k^2$, thus implying that $f(x_k)$ will finally oscillate around the origin, which is a contradiction with the linear convergence of $f(x_k)$. Hence, (\ref{10}) is not linearly convergent.
\end{example}

\begin{example}\label{e3}
Consider now another sign variant as
\begin{align}
b_{k+1}^2&=b_k^2+\|\operatorname{sign}(\nabla f(x_k))\|^2,         \nonumber\\
x_{k+1}&=x_k-\frac{\eta}{b_{k+1}}\operatorname{sign}(\nabla f(x_k)),~~~\eta>0        \label{14}
\end{align}
and let $f=x^2$ (strongly convex and smooth) with $x\in\mathbb{R}$ with initial $x_0=\eta/2$ and $b_0=0$. In this case, it is straightforward to calculate that $b_k=\sqrt{k}$ and
\begin{align}
x_{k}&=\eta\Big(\frac{1}{2}-1+\frac{1}{\sqrt{2}}-\frac{1}{\sqrt{3}}+\cdots+\frac{(-1)^{k}}{\sqrt{k}}\Big),     \label{15}
\end{align}
from which one can conclude that (\ref{10}) amounts to
\begin{align}
x_{k+1}=x_k-\frac{\eta}{\sqrt{k+1}}\operatorname{sign}(x_k),        \label{16}
\end{align}
which can be viewed as GD for the convex objective $g(x)=|x|$ with a learning rate $\eta/\sqrt{k+1}$. Therefore, the convergence rate of classic GD can be invoked for (\ref{13}), which is known to be sublinear \cite{nedic2015distributed}.
\end{example}

\begin{remark}\label{r1}
The above examples demonstrate that although GD and AdaGrad-Norm are indeed linearly convergent for strongly convex and smooth objectives, their sign counterparts fail to converge linearly in general.
\end{remark}

\section{Linear Rate of Scaled {\footnotesize SIGN}GD/SGD}\label{s4}

With the above preparations, it is now ready to study faster convergence for solving problem (\ref{1}). As shown in Section \ref{s3}, the sign counterparts of GD and AdaGrad-Norm are not applicable for linear convergence. As such, the scaled versions of {\scriptsize SIGN}GD/SGD are considered in this paper, as in Algorithms \ref{al1} and \ref{al2}, which can be viewed as the steepest descent with respect to the maximum norm \cite{balles2020geometry}, but is still not fully understood.

A few assumptions are necessary for the following analysis.

\begin{assumption}\label{a1}
$f$ is $\mu$-strongly convex with respect to $\ell_\infty$-norm for some constant $\mu>0$, i.e., $f(x)-f(y)\geq \nabla f(y)^\top(x-y)+\frac{\mu}{2}\|x-y\|_\infty^2$ for all $x,y\in\mathbb{R}^d$.
\end{assumption}
\begin{assumption}\label{a2}
$f$ satisfies the Polyak-{\L}ojasiewicz (PL) inequality, i.e., $\|\nabla f(x)\|_1^2\geq 2\mu(f(x)-f^*),~\forall x\in\mathbb{R}^d$, where $f^*$ is the minimum value.
\end{assumption}
\begin{assumption}\label{a3}
$f$ is $L$-smooth with respect to $\ell_\infty$-norm, i.e., $\|\nabla f(x)-\nabla f(y)\|_1\leq L\|x-y\|_\infty$ for all $x,y\in\mathbb{R}^d$.
\end{assumption}

The PL inequality does not require $f$ to be even convex, and the $\ell_\infty$- and $\ell_1$-norms employed in Assumptions \ref{a1} and \ref{a2}, respectively, are slightly more relaxed than the Euclidean norm. Meanwhile, the smoothness condition is made with respect to $\ell_\infty$-norm, since it is more favorable than the Euclidean smoothness and separable smoothness \cite{balles2020geometry}.

\begin{remark}\label{r2}
It is noteworthy that another promising sign method is {\scriptsize EF-SIGN}GD \cite{karimireddy2019error} using error feedback, given as
\begin{align}
p_k&=\lambda\nabla f(x_k)+e_k,        \nonumber\\
x_{k+1}&=x_k-\frac{\|p_k\|_1}{d}\operatorname{sign}(p_k),           \nonumber\\
e_{k+1}&=p_k-\frac{\|p_k\|_1}{d}\operatorname{sign}(p_k),           \label{17}
\end{align}
where $\lambda>0$ is the learning rate. In \cite{karimireddy2019error}, it is shown that {\scriptsize EF-SIGN}GD/SGD has a better performance than {\scriptsize SIGN}GD/SGD, actually enjoying the same convergence rate as GD/SGD. However, we point out that {\scriptsize EF-SIGN}GD/SGD is, roughly speaking, equivalent to GD/SGD. Let us show this by slightly modifying (\ref{17}) as
\begin{align}
p_k&=\lambda\nabla f(x_k-e_k)+e_k,        \nonumber\\
x_{k+1}&=x_k-\frac{\|p_k\|_1}{d}\operatorname{sign}(p_k),           \nonumber\\
e_{k+1}&=p_k-\frac{\|p_k\|_1}{d}\operatorname{sign}(p_k).           \label{18}
\end{align}
By defining $z_k=x_k-e_k$, it is easy to verify that $z_{k+1}=z_k-\lambda\nabla f(z_k)$, that is, (\ref{18}) amounts to GD in terms of $z_k$. As a result, {\scriptsize EF-SIGN}GD/SGD is not considered here.
\end{remark}

In the following, the main results are divided into two scenarios, i.e., the deterministic and stochastic settings.

\begin{algorithm}[tb]
 \caption{Scaled {\scriptsize SIGN}GD}   \label{al1}
 \begin{algorithmic}
  \STATE \textbf{Input:} learning rate $\alpha$, current point $x_k$
  \STATE \vspace{-0.6cm}
  \begin{align}
x_{k+1}=x_k-\alpha \|g_k\|_1 \operatorname{sign}(g_k),~~~g_k:=\nabla f(x_k)           \label{4}
\end{align}
 \end{algorithmic}
\end{algorithm}

\begin{algorithm}[tb]
 \caption{Scaled {\scriptsize SIGN}SGD}      \label{al2}
 \begin{algorithmic}
  \STATE \textbf{Input:} learning rate $\alpha_k$, current point $x_k$
  \STATE \vspace{-0.6cm}
\begin{align}
&\hspace{-2.7cm}\tilde{g}_k=\text{StochasticGradient}(x_k)            \nonumber\\
&\hspace{-2.7cm}x_{k+1}=x_k-\alpha_k \|\tilde{g}_k\|_1\operatorname{sign}(\tilde{g}_k)           \label{5}
\end{align}
 \end{algorithmic}
\end{algorithm}

\subsection{The Deterministic Setting}\label{s4.1}

Consider the deterministic setting with full gradients, i.e., (\ref{4}), for which we have the following results. Note that all proofs are given in the Appendix.

\begin{theorem}\label{t1}
The following statements are true for (\ref{4}).
\begin{enumerate}
  \item Under Assumptions \ref{a1} and \ref{a3}, if $0<\alpha<\frac{2}{L}$, then
  \begin{align}
  f(x_k)-f^*\leq \zeta^k(f(x_0)-f^*),         \label{20}
  \end{align}
  where $\zeta:=1-2\mu\alpha\big(1-\frac{L\alpha}{2}\big)\in[0,1)$.
  \item Under Assumptions \ref{a2} and \ref{a3} with $\alpha$ satisfying $0<\alpha<\frac{2}{L}$, (\ref{20}) still holds.
  \item If Assumption \ref{a3} holds only, then
  \begin{align}
  \min_{l\in\{0,1,\ldots,k\}}\|g_l\|_1^2\leq \frac{f(x_0)-f^*}{\gamma (k+1)},         \label{21}
  \end{align}
  where $\gamma:=\alpha(1-\frac{L\alpha}{2})$.
\end{enumerate}
\end{theorem}

\begin{remark}\label{r3}
In view of Theorem \ref{t1}, the algorithm (\ref{4}) is proved to be linearly convergent, which is contrast to {\scriptsize SIGN}GD and sign AdaGrad-Norm as discussed in Section \ref{s3}. Moreover, for the nonconvex but smooth with respect to the Euclidean norm, by leveraging the similar argument to Theorem \ref{t1}, it is easy to obtain for {\scriptsize SIGN}GD with a constant learning rate that $\min_{l\in\{0,1,\ldots,k\}}\|g_l\|_1^2\leq \frac{dL(f(x_0)-f^*)}{2(k+1)}$ by choosing the learning rate as $\alpha=\sqrt{\frac{2(f(x_0)-f^*)}{dL(k+1)}}$. In comparison, (\ref{21}) can be nearly $\frac{L(f(x_0)-f^*)}{2(k+1)}$ when $\alpha$ is chosen to approach $\frac{2}{L}$. In this regard, our result is tighter up to a dimension constant $d$, and the learning rate here is easier to implement. In addition, if the smoothness is with respect to the maximum norm, then the result here has the same convergence bound as {\scriptsize SIGN}GD but with a less conservative learning rate selection.
\end{remark}

\begin{remark}\label{rr1}
A similar result can be also obtained from the most related work \cite{beznosikov2020biased} by resorting to the $\delta$-approximate compressor. To be specific, $\mathcal{C}(v):=\frac{\|v\|_1}{d}\operatorname{sign}(v)$ can be viewed as $\frac{1}{d}$-approximate compressor, and then applying Theorem 13 in \cite{beznosikov2020biased} leads to the learning rate $\alpha\in[0,\frac{1}{L}]$ and convergence rate $(1-\frac{\alpha\mu}{d})^k$. In contrast, Theorem 1 of this paper (need to replace $\alpha$ by $\frac{\alpha}{d}$ here) is for $\alpha\in(0,\frac{2d}{L})$ with the convergence rate $(1-\frac{2\mu\alpha}{d}(1-\frac{L\alpha}{2}))^k$. It is easy to verify that our learning rate is more relaxed and the convergence rate is faster due to $\frac{\alpha\mu}{d}\leq \frac{2\mu\alpha}{d}(1-\frac{L\alpha}{2})$.
\end{remark}

\subsection{The Stochastic Setting}\label{s4.2}

This section considers the stochastic gradient case, where the true gradient $g_k=\nabla f(x_k)$ is expensive to compute and instead a stochastic gradient $\tilde{g}_k$ is relatively cheap to evaluate as an estimate of $g_k$. To move forward, some standard assumptions are imposed on stochastic gradients \cite{balles2018dissecting,bernstein2018signsgd}.
\begin{assumption}\label{a4}
The stochastic gradients $\{\tilde{g}_k\}_{k=0}^\infty$ are unbiased and have bounded variances with respect to $\ell_1$-norm, i.e., there exists a constant $\sigma>0$ such that
\begin{align}
\mathbb{E}(\tilde{g}_k)=g_k,~~~~~\mathbb{E}(\|\tilde{g}_k-g_k\|_1^2)\leq \sigma^2.            \label{24}
\end{align}
\end{assumption}

In this case, the algorithm becomes (\ref{5}). For brevity, define $p_{k,i}:=\mathbb{P}(\operatorname{sign}(\tilde{g}_{k,i})=\operatorname{sign}(g_{k,i}))$ for $k\geq 0$ and $i\in[d]$, where $\tilde{g}_{k,i}$ and $g_{k,i}$ represents the $i$-th components of $\tilde{g}_{k,i}$ and $g_k$, respectively.

\begin{remark}\label{r4}
For stochastic gradient $g_k$, when leveraging a mini-batch of size $n_k$ at $x_k$, the oracle gives us $n_k$ gradient estimates and in this case, the stochastic gradient $\tilde{g}_k$ can be chosen as the average of $n_k$ estimates. In this respect, the variance bound can be reduced to $\frac{\sigma^2}{n_k}$. Additionally, it was shown in \cite{safaryan2019stochastic} that the success probability $p_{k,i}$ should be greater than $1/2$, and otherwise the sign algorithm generally fails to work. And a multitude of cases can ensure $p_{k,i}>1/2$, for instance, each component $\tilde{g}_{k,i}$ possesses a unimodal and symmetric distribution \cite{bernstein2018signsgd,safaryan2019stochastic}.
\end{remark}

We are now in a position to present the main result on (\ref{5}).
\begin{theorem}\label{t3}
For (\ref{5}), under Assumptions \ref{a1}, \ref{a3}, \ref{a4} or \ref{a2}-\ref{a4}, the following statements are true.
\begin{enumerate}
  \item If $\alpha_k=\alpha\in\big(0,\frac{2p_{min}-1}{L}\big)$, then
\begin{align}
\mathbb{E}(f(x_k))-f^*&\leq \zeta_1^k(\mathbb{E}(f(x_0))-f^*)             \nonumber\\
&\hspace{0.4cm}+\frac{L\sigma^2\alpha}{2\mu(2p_{min}-1-L\alpha)},       \label{26}
\end{align}
where $p_{min}:=\min_{i\in[d],k\geq 0}p_{k,i}$ and $\zeta_1:=1-2\mu\alpha(2p_{min}-1-L\alpha)\in[\frac{1}{2},1)$.
  \item If $\alpha_k=\frac{3}{\mu(2p_{min}-1)(k+1)}$, then
  \begin{align}
\mathbb{E}(f(x_k))-f^*&\leq \frac{9L\sigma^2}{\mu^2(2p_{min}-1)^2}\Big(\frac{32}{k}+\frac{1}{k^2}\Big)    \nonumber\\
&\hspace{0.4cm}+\frac{f(x_0)-f^*}{(k+1)^3}.       \label{27}
  \end{align}
\end{enumerate}
\end{theorem}

\begin{remark}\label{r5}
The first result in Theorem \ref{t3} shows that algorithm (\ref{5}) converges linearly at a rate $\zeta_1$. This is comparable to vanilla SGD in \cite{gower2019sgd}, where the convergence rate is $1-\alpha\mu$, which is slower than $\zeta_1$ (i.e., $\zeta_1\leq 1-\alpha\mu$) when $\alpha\in(0,\frac{1}{2L})$. Moreover, the result in (\ref{27}) is the exact convergence with rate $O(\frac{1}{k})$ for both strongly convex case and nonconvex case with PL inequality, which is the same as both vanilla SGD \cite{rakhlin2012making} and compression methods \cite{stich2018sparsified}. In addition, the same rate $O(\frac{1}{k})$ was established in \cite{carlson2015stochastic}. However, the condition in \cite{carlson2015stochastic} for convergence does not always hold, e.g., $t_k=1$ in Theorem II.2 of \cite{carlson2015stochastic}, and our result (\ref{27}) includes more faster rate $O(\frac{1}{k^2}+\frac{1}{k^3})$ except for $O(\frac{1}{k})$ in \cite{carlson2015stochastic}.
\end{remark}


\section{The Distributed Setting}\label{s5}

Now, we extend the results in Section \ref{s4} to the distributed setting within a parameter server framework. For simplicity, we only focuses on scaled {\scriptsize SIGN}SGD in this section, but the results can be similarly obtained for scaled {\scriptsize SIGN}GD.

\begin{algorithm}[tb]
 \caption{Distributed Scaled {\scriptsize SIGN}SGD by Majority Vote}
 \label{al3}
 \begin{algorithmic}
  \STATE \textbf{Input:} learning rate $\alpha$, current point $x_k$, $\#$ workers $M$ each with an i.i.d. gradient estimate $\tilde{g}_k^m,m\in[M]$
  \STATE \textbf{On} server

  \hspace{0.8cm}{\bf Pull} $\operatorname{sign}(\tilde{g}_k^m)$ and $\|\tilde{g}_k^m\|_1$ {\bf from} each worker\\
  \hspace{0.8cm}{\bf Push} $\operatorname{sign}(\hat{g}_k^s)$ and $M_k$ {\bf to} each worker\\
  \hspace{1.7cm}$\hat{g}_k^s:=\frac{1}{M}\sum_{m=1}^M\operatorname{sign}(\tilde{g}_k^m)$\\
  \hspace{1.7cm}$M_k:=\frac{1}{M}\sum_{m=1}^M\|\tilde{g}_k^m\|_1$
  \STATE \textbf{On} each worker

  \hspace{0.8cm}$x_{k+1}=x_k-\alpha M_k\operatorname{sign}(\hat{g}_k^s)$
 \end{algorithmic}
\end{algorithm}

To proceed, the distributed scaled {\scriptsize SIGN}SGD by majority vote is given in Algorithm \ref{al3}, for which the following convergence result is obtained.

\begin{theorem}\label{t4}
For Algorithm \ref{al3}, under Assumptions \ref{a1}, \ref{a3}, \ref{a4} or \ref{a2}-\ref{a4}, if $0<\alpha <\frac{2I_{p_{min}}(\kappa,\kappa)-1}{L}$, then
\begin{align}
\mathbb{E}(f(x_k))-f^*\leq \zeta_2^k(\mathbb{E}(f(x_0))-f^*)+\frac{L\sigma^2\alpha^2}{1-\zeta_2},       \label{29}
\end{align}
where $p_{min}=\min_{l\in[d],k\geq 0}p_{k,l}$, $\zeta_2:=1-2\mu\alpha(2I_{p_{min}}(\kappa,\kappa)-1-L\alpha)\in[\frac{1}{2},1)$, $\kappa:=\lfloor\frac{M+1}{2}\rfloor$ with $\lfloor\cdot\rfloor$ being the floor function, and $I_p(a,b)$ is the regularized incomplete beta function, defined by
\begin{align}
I_p(a,b):=\frac{\int_0^p t^{a-1}(1-t)^{b-1}dt}{\int_0^1 t^{a-1}(1-t)^{b-1}dt},~~~a,b>0,~p\in[0,1].        \nonumber
\end{align}
\end{theorem}

\begin{remark}\label{r6}
It is noteworthy that the exact convergence can be similarly established as (\ref{27}) in Theorem \ref{t3}, which is omitted in Theorem \ref{t4}.
\end{remark}

\section{Experiments}

Numerical experiments are provided to corroborate the efficacy of the obtained theoretical results here.

\begin{figure}[H]
\centering
\includegraphics[width=2.6in]{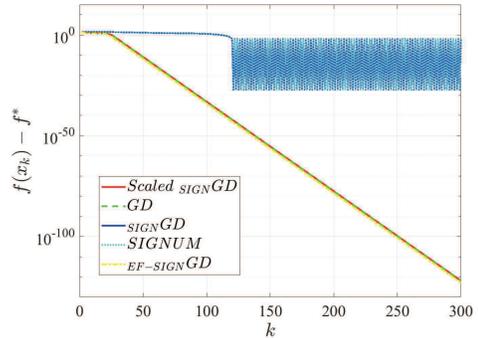}
\caption{Simulation results for several algorithms.}
\label{f0}
\end{figure}

\begin{example}[A Toy Example]
Let us consider a simple example, where $f(x)=x^2+3\sin^2(x)$ for $x\in\mathbb{R}$. It is easy to verify that $f(x)$ is nonconvex, but satisfying the PL condition. To verify the performance of the proposed scaled {\scriptsize SIGN}GD, several existing algorithms are compared in Fig. \ref{f0} by setting $\alpha=0.05$ with an arbitrary initial state. The comparisons are performed with vanilla gradient descent (GD), {\scriptsize SIGN}GD, {\scriptsize SIGN}GD{\scriptsize M} (i.e., SIGNUM), and {\scriptsize EF-SIGN}GD \cite{karimireddy2019error}. It can be observed from Fig. \ref{f0} that the proposed algorithm has the same linear convergence as GD and {\scriptsize EF-SIGN}GD, while {\scriptsize SIGN}GD and SIGNUM cannot converge, behaving oscillations near the optimal variable. In summary, this example shows the efficiency of the scaled {\scriptsize SIGN}GD, and supports the observation in Example \ref{e1}.
\end{example}

\begin{example}
Consider the logistic regression problem, where the objective is $f(x)=\frac{1}{n}\sum_{i=1}^n \log(1+\exp(-b_ia_i^\top x))+\frac{1}{2n}\|x\|^2$ with a standard $L2$-regularizer \cite{stich2018sparsified}, and $a_i\in\mathbb{R}^d$ and $b_i\in\{-1,+1\}$ are the data samples.

\begin{figure}[H]
\centering
\subfigure[]{\includegraphics[width=2.5in]{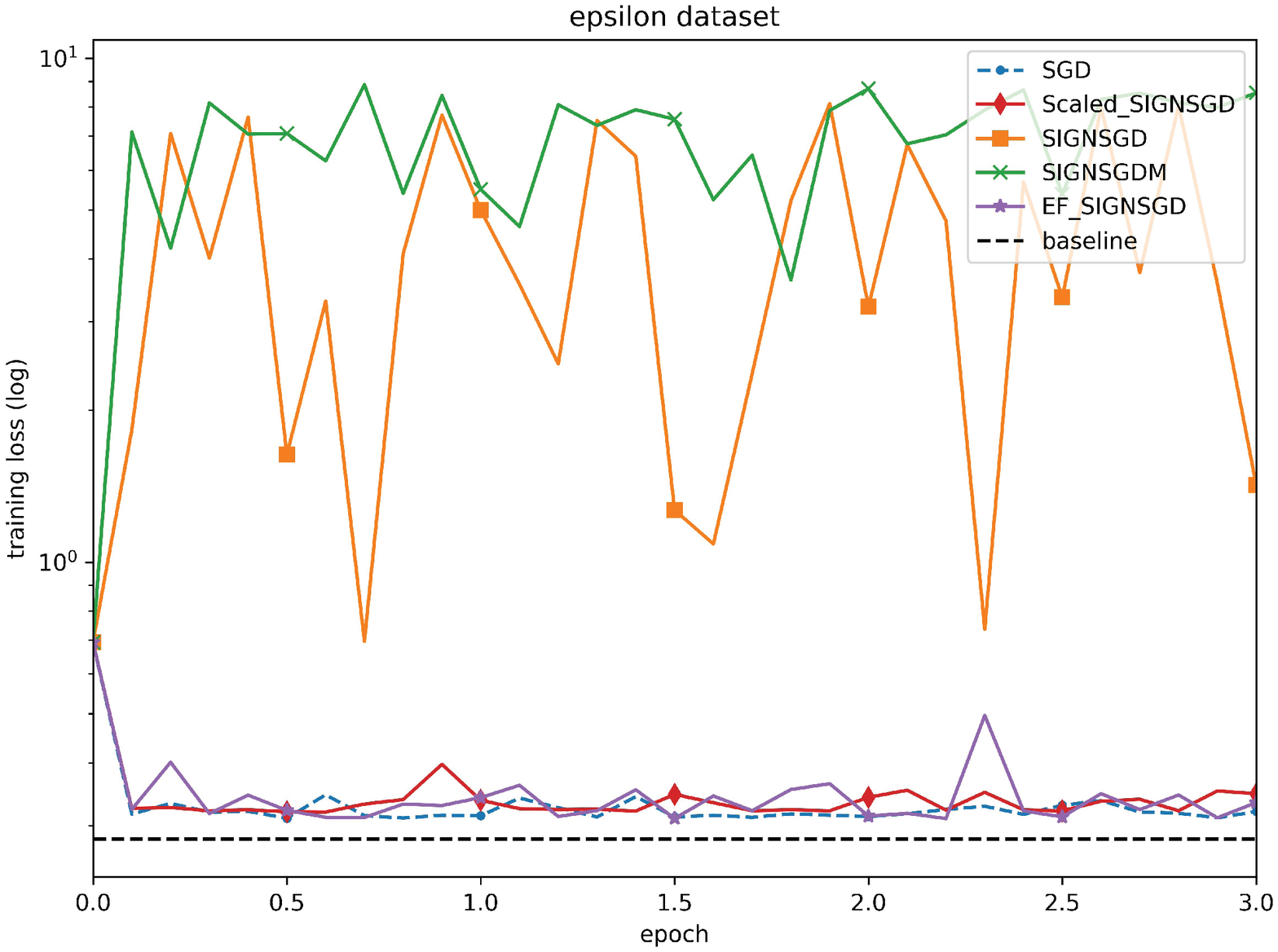}}\hspace{0.0cm}
\subfigure[]{\includegraphics[width=2.5in]{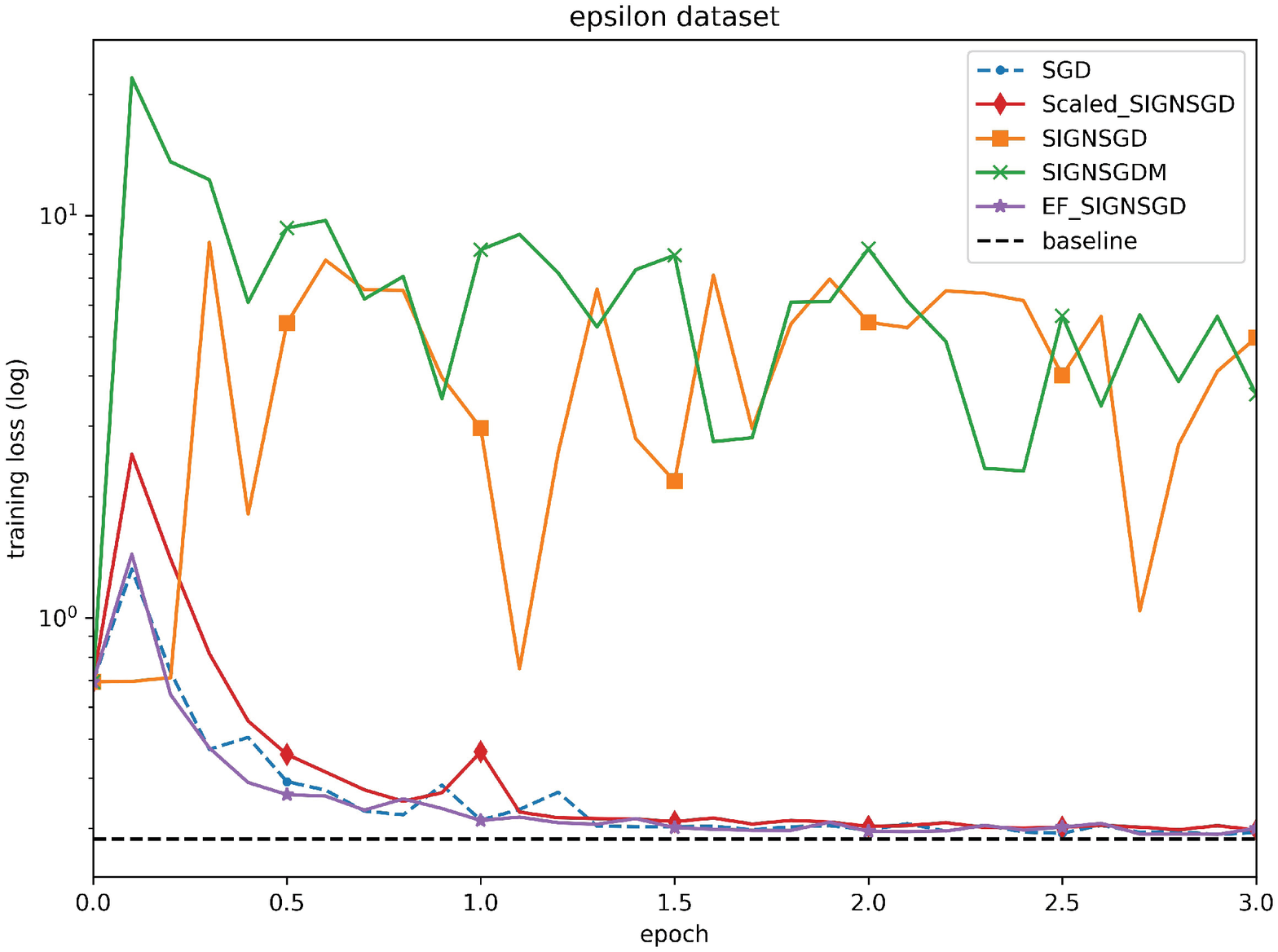}}
\caption{Scaled {\scriptsize SIGN}SGD. (a) $\alpha=2$; (b) $\alpha_k=\frac{6n}{k+1}$.}
\label{f1}
\end{figure}

\begin{figure}[H]
\centering
\includegraphics[width=2.6in]{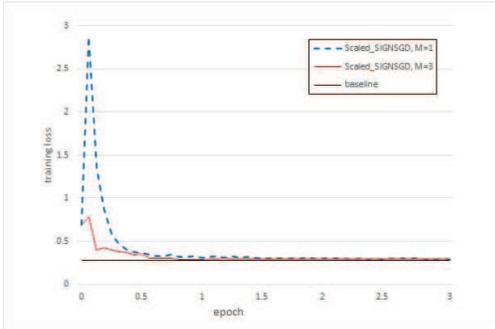}
\caption{Distributed scaled {\scriptsize SIGN}SGD for $M=1,3$.}
\label{f2}
\end{figure}

To test the performance of scaled {\scriptsize SIGN}SGD, the {\em epsilon} dataset with $n=400000$ and $d=2000$ is exploited \cite{sonnenburg2008pascal}, and the baseline is calculated using the standard optimizer LogisticSGD of scikit-learn \cite{pedregosa2011scikit}. To marginalize out the effect of initial choices, the numerical result is averaged over repeated runs with $x_0\approx \mathcal{N}(0,I)$. We compare scaled {\scriptsize SIGN}SGD with vanilla SGD, {\scriptsize SIGN}SGD, {\scriptsize SIGN}SGD{\scriptsize M}, and {\scriptsize EF-SIGN}SGD \cite{karimireddy2019error}, as shown in Fig. \ref{f1} on a platform with the Intel Core i7-4300U CPU. Fig. \ref{f1} indicates that {\scriptsize SIGN}SGD has a similar performance to SGD and performs better than {\scriptsize SIGN}SGD and {\scriptsize SIGN}SGD{\scriptsize M}. It can be also observed that {\scriptsize EF-SIGN}SGD is comparable to SGD, which is consistent with the discussion in Remark \ref{r2}. Moreover, the case in Fig. \ref{f1}(a) with a constant learning rate converges faster than that in Fig. \ref{f1}(b) with a diminishing learning rate. Meanwhile, Fig. \ref{f2} shows that more workers can improve the performance. Therefore, the numerical results support our theoretical findings.
\end{example}

\section{Conclusion}

This paper has investigated faster convergence of scaled {\scriptsize SIGN}GD/SGD, which can relieve the communication cost compared with vanilla SGD. To further motivate the study of sign methods, continuous-time algorithms have been addressed, indicating that sign SGD can significantly improve the convergence speed of SGD. Subsequently, it has been proven that scaled {\scriptsize SIGN}GD is linearly convergent for both strongly convex and nonconvex (satisfying PL inequality) objectives. Also, the convergence for {\scriptsize SIGN}SGD has been analyzed in two cases with constant and decaying learning rates. The results are also extended to the distributed setting in the parameter server framework. The efficacy of scaled sign methods has been validated by numerical experiments for the logistic regression problem.

%

\section*{Appendix}

\subsection{Further Motivations for {\scriptsize SIGN}GD}\label{s2}

Let us provide more evidences for studying sign-based GD from the continuous-time perspective. In doing so, consider the continuous-time dynamics corresponding to the discrete-time GD and {\scriptsize SIGN}GD, i.e.,
\begin{align}
\dot{x}&=-\beta \nabla f(x),          \label{ca1}\\
\dot{x}&=-\beta \operatorname{sign}(\nabla f(x)),        \label{ca2}
\end{align}
where $\beta>0$ is a constant learning rate.

To proceed, let us construct a Lyapunov candidate as
\begin{align}
V(t):=f(x)-f^*,~~~\forall~t\geq 0       \label{6}
\end{align}
where $f^*$ denotes the minimum value attained by $f$.

For algorithms (\ref{ca1}) and (\ref{ca2}), the following results can be obtained.

\begin{proposition}\label{p1}
For algorithm (\ref{ca1}),
\begin{enumerate}
  \item if $f$ is convex, then $V(t)\leq \frac{D_1^2 V(0)}{D_1^2+V(0)\beta t}$, where $D_1:=\max_{x:f(x)\leq f(x_0)}\min_{x^*\in\mathcal{X}^*}\|x-x^*\|$ with $\mathcal{X}^*$ being the set of minimizers;
  \item if $f$ is nonconvex, then $\min_{s\in[0,t]}\|\nabla f(x(s))\|\leq \frac{\sqrt{f(x(0))-f^*}}{\sqrt{\beta t}}$.
\end{enumerate}
\end{proposition}

\begin{proposition}\label{p2}
For algorithm (\ref{ca2}),
\begin{enumerate}
  \item if $f$ is convex, then $V(t)\leq V(0)e^{-\frac{\beta t}{D_2}}$, where $D_2:=\max_{x:f(x)\leq f(x_0)}\min_{x^*\in\mathcal{X}^*}\|x-x^*\|_\infty$;
  \item if $f$ is nonconvex, then $\min_{s\in[0,t]}\|\nabla f(x(s))\|_1\leq \frac{f(x(0))-f^*}{\beta t}$.
\end{enumerate}
\end{proposition}

In view of the above results, it can be easily observed that (\ref{ca2}) with sign gradients converges apparently faster than GD (\ref{ca1}) in the continuous-time domain, indicating that the performance of gradient descent can be largely improved by sign gradient compression. For instance, in the scenario with convex objectives, GD (\ref{ca1}) is sublinearly convergent while {\scriptsize SIGN}GD (\ref{ca2}) is linearly convergent. As a result, the above results provide a new perspective for showing advantages of {\scriptsize SIGN}GD compared with GD.

\subsection{Proof of Proposition \ref{p1}}

Consider the case with convex objectives. In light of (\ref{ca1}), it can be calculated that
\begin{align}
\dot{V}(t)=\nabla f(x(t))^\top\dot{x}=-\beta\|\nabla f(x(t))\|^2\leq 0,         \label{pf1}
\end{align}
which implies $f(x(t))\leq f(x(0))$.

Meanwhile, invoking the convexity of $f$ yields
\begin{align}
V(t)&\leq\nabla f(x(t))^\top (x-x^*)        \nonumber\\
&\leq \|\nabla f(x(t))\|\cdot\|x-x^*\|       \nonumber\\
&\leq D_1\|\nabla f(x(t))\|,           \label{pf3}
\end{align}
which, together with (\ref{pf1}), gives rise to $\dot{V}(t)\leq -\frac{\beta}{D_1^2}V(t)^2$, further implying the claimed result.

For the case with nonconvex objectives, by integrating (\ref{pf1}) from $0$ to $t$, one can obtain that
\begin{align}
\beta\int_0^t\|\nabla f(x(s))\|^2ds&=V(0)-V(t)         \nonumber\\
&=f(x(0))-f(x(t))            \nonumber\\
&\leq f(x(0))-f^*,           \label{pf5}
\end{align}
where the inequality has employed the fact that $f(z)\geq f^*$ for all $z\in\mathbb{R}^d$. Then taking the minimum of $\|\nabla f(x(s))\|$ over $[0,t]$ ends the proof.
\qed

\subsection{Proof of Proposition \ref{p2}}

Consider first the convex case. Similar to (\ref{pf1}), it can be obtained that
\begin{align}
\dot{V}(t)=-\beta\|\nabla f(x(t))\|_1.            \label{pf6}
\end{align}
Akin to (\ref{pf3}), one has that
\begin{align}
V(t)&\leq\nabla f(x(t))^\top (x-x^*)        \nonumber\\
&\leq \|\nabla f(x(t))\|_1\cdot\|x-x^*\|_\infty       \nonumber\\
&\leq D_2\|\nabla f(x(t))\|_1,           \label{pf7}
\end{align}
where the second inequality has used Holder's inequality. Combining (\ref{pf6}) with (\ref{pf7}) yields $\dot{V}(t)\leq -\frac{\beta}{D_2}V(t)$, from which it is easy to verify the claimed result.

Consider now the nonconvex case. The desired result can be obtained by (\ref{pf6}) and the similar argument to that in convex case. This completes the proof.
\qed

\subsection{Proof of Theorem \ref{t1}}

To facilitate the subsequent analysis, define
\begin{align}
V_k:=f(x_k)-f^*,~~~\forall k\geq 0.        \label{pf9}
\end{align}
In view of (\ref{4}) and Assumption \ref{a3}, it can be concluded that
\begin{align}
V_{k+1}-V_k&=f(x_{k+1})-f(x_k)        \nonumber\\
&\leq \nabla f(x_k)^\top(x_{k+1}-x_k)+\frac{L}{2}\|x_{k+1}-x_k\|_\infty^2        \nonumber\\
&=-\alpha \|g_k\|_1^2+\frac{L\alpha^2}{2}\|g_k\|_1^2\cdot\|\operatorname{sign}(g_k)\|_\infty^2     \nonumber\\
&\leq -\alpha \|g_k\|_1^2+\frac{L\alpha^2}{2}\|g_k\|_1^2            \nonumber\\
&=-\gamma \|g_k\|_1^2.                  \label{pf10}
\end{align}
In what follows, let us prove this theorem one by one.

First, for case 1, invoking Assumption \ref{a1} yields
\begin{align}
V_k&\leq g_k^\top(x_k-x^*)-\frac{\mu}{2}\|x_k-x^*\|_\infty^2         \nonumber\\
&\leq \frac{1}{2}\Big(\frac{\|g_k\|_1^2}{\mu}+\mu\|x_k-x^*\|_\infty^2\Big)-\frac{\mu}{2}\|x_k-x^*\|_\infty^2          \nonumber\\
&=\frac{\|g_k\|_1^2}{2\mu},          \nonumber
\end{align}
where the second inequality has employed the Holder inequality. Then one has that $\|g_k\|_1^2\geq 2\mu V_k$. Therefore, in combination with (\ref{pf10}), one can obtain that $V_{k+1}-V_k\leq -2\mu\gamma V_k$, further leading to $V_{k+1}\leq \zeta V_k$. Consequently, by iteration, this completes the proof of case 1.

Second, for case 2, Assumption \ref{a2} leads to $2\mu V_k\leq \|g_k\|_1^2$, which, together with the similar argument to case 1, follows the conclusion in this case.

Third, for case 3, invoking (\ref{pf10}) gives $\gamma \|g_k\|_1^2\leq V_k-V_{k+1}$, which, by summation over $l=0,1,\ldots,k$, implies that
\begin{align}
\gamma\sum_{l=0}^{k}\|g_l\|_1^2&\leq V_0-V_{k+1}=f(x_0)-f(x_{k+1})                \nonumber\\
&\leq f(x_0)-f^*,              \label{pf11}
\end{align}
where the last inequality has used the fact that $f(x_{k+1})\geq f^*$. Then taking the minimum of $\|g_l\|_1^2$ over $l=0,1,\ldots,k$ ends the proof.
\qed

\subsection{Proof of Theorem \ref{t3}}

Recalling $V_k$ in (\ref{pf9}). Invoking Assumption \ref{a3} gives rise to
\begin{align}
V_{k+1}-V_k&\leq g_k^\top(x_{k+1}-x_k)+\frac{L}{2}\|x_{k+1}-x_k\|_\infty^2          \nonumber\\
&\hspace{-0.9cm}=-\alpha\|\tilde{g}_k\|_1g_k^\top\operatorname{sign}(\tilde{g}_k)+\frac{L\alpha^2}{2}\|\tilde{g}_k\|_1^2\|\operatorname{sign}(\tilde{g}_k)\|_\infty^2       \nonumber\\
&\hspace{-0.9cm}\leq -\alpha\|\tilde{g}_k\|_1g_k^\top\operatorname{sign}(\tilde{g}_k)+\frac{L\alpha^2}{2}\|\tilde{g}_k\|_1^2.         \nonumber
\end{align}
By taking the conditional expectation, one has
\begin{align}
\mathbb{E}(V_{k+1}|x_k)-V_k&\leq -\alpha g_k^\top\mathbb{E}(\|\tilde{g}_k\|_1\operatorname{sign}(\tilde{g}_k)|x_k)            \nonumber\\
&\hspace{0.4cm}+\frac{L\alpha^2}{2}\mathbb{E}(\|\tilde{g}_k\|_1^2|x_k).         \label{pf14}
\end{align}

Consider now the coordinate $\tilde{g}_{k,i}$ for $i\in[d]$. One has that
\begin{align}
\mathbb{E}(\|\tilde{g}_k\|_1\operatorname{sign}(\tilde{g}_{k,i})|x_k)&=\mathbb{E}[\mathbb{E}(\|\tilde{g}_k\|_1\operatorname{sign}(\tilde{g}_{k,i})|\tilde{g}_k)|x_k]     \nonumber\\
&\hspace{-3.4cm}=\mathbb{E}[\|\tilde{g}_k\|_1\mathbb{E}(\operatorname{sign}(\tilde{g}_{k,i})|\tilde{g}_k)|x_k]           \nonumber\\
&\hspace{-3.4cm}=\mathbb{E}[\|\tilde{g}_k\|_1\mathbb{P}(\operatorname{sign}(\tilde{g}_{k,i})=\operatorname{sign}(g_{k,i}))\operatorname{sign}(g_{k,i})          \nonumber\\
&\hspace{-3.0cm}-\|\tilde{g}_k\|_1\mathbb{P}(\operatorname{sign}(\tilde{g}_{k,i})\neq\operatorname{sign}(g_{k,i}))\operatorname{sign}(g_{k,i})|x_k]            \nonumber\\
&\hspace{-3.4cm}=\mathbb{E}[p_{k,i}\|\tilde{g}_k\|_1\operatorname{sign}(g_{k,i})-(1-p_{k,i})\|\tilde{g}_k\|_1\operatorname{sign}(g_{k,i})|x_k]                   \nonumber\\
&\hspace{-3.4cm}=(2p_{k,i}-1)\operatorname{sign}(g_{k,i})\mathbb{E}(\|\tilde{g}_k\|_1|x_k),         \nonumber
\end{align}
which, together with (\ref{pf14}), implies that
\begin{align}
\mathbb{E}(V_{k+1}|x_k)-V_k&\leq -\alpha \sum_{i=1}^d (2p_{k,i}-1)|g_{k,i}|\mathbb{E}(\|\tilde{g}_k\|_1|x_k)            \nonumber\\
&\hspace{0.4cm}+\frac{L\alpha^2}{2}\mathbb{E}(\|\tilde{g}_k\|_1^2|x_k).                \label{pf15}
\end{align}
By Jesen's inequality, it follows that $\mathbb{E}(\|\tilde{g}_k\|_1|x_k)\geq \|\mathbb{E}(\tilde{g}_k|x_k)\|_1=\|g_k\|_1$. Because $p_{k,i}\geq p_{min}$ for $i\in[d]$, taking the expectation implies that
\begin{align}
\mathbb{E}(V_{k+1})-\mathbb{E}(V_k)&\leq -\alpha (2p_{min}-1)\mathbb{E}(\|g_k\|_1^2)            \nonumber\\
&\hspace{-2.0cm}+L\alpha^2[\mathbb{E}(\|\tilde{g}_k-g_k\|_1^2)+\mathbb{E}(\|g_k\|^2)]        \nonumber\\
&\hspace{-2.4cm}\leq -\alpha (2p_{min}-1-L\alpha)\mathbb{E}(\|g_k\|_1^2)+L\sigma^2\alpha^2.                 \label{pf16}
\end{align}
Now, under Assumption \ref{a1} or \ref{a2}, using the similar argument to the proof of Theorem \ref{t1} can both lead to that $\mathbb{E}(\|g_k\|_1^2)\geq 2\mu\mathbb{E}(V_k)$, which together with (\ref{pf16}) yields that
\begin{align}
\mathbb{E}(V_{k+1})\leq \zeta_1\mathbb{E}(V_k)+L\sigma^2\alpha^2.          \label{new1}
\end{align}
Iteratively applying the above inequality leads to (\ref{26}).

It remains to show (\ref{27}). Invoking the similar analysis for (\ref{new1}) yields that
\begin{align}
\mathbb{E}(V_{k+1})\leq c_k\mathbb{E}(V_k)+L\sigma^2\alpha_k^2,           \nonumber
\end{align}
where $c_k:=1-\alpha_k \mu(2p_{min}-1)$, further implying that
\begin{align}
\mathbb{E}(V_k)&\leq \Pi_{l=0}^{k-1} c_l V_0+L\sigma^2(c_{k-1}\cdots c_1\alpha_0^2+\cdots       \nonumber\\
&\hspace{0.4cm}+c_{k-1}\alpha_{k-2}^2+\alpha_{k-1}^2)            \nonumber\\
&\leq \Pi_{l=0}^{k-1} c_l V_0+\frac{9}{\mu^2(2p_{min}-1)^2}          \nonumber\\
&\hspace{1.5cm}\cdot\Big(\sum_{m=1}^{k-1}e^{-\sum_{l=m}^{k-1}\frac{3}{l+1}}\frac{1}{m^2}+\frac{1}{k^2}\Big),           \label{new2}
\end{align}
where the second inequality has employed the expression of $\alpha_k$.

For the last two terms in (\ref{new2}), in light of the fact that $\Pi_{l=0}^m(1-a_l)\leq e^{-\sum_{l=0}^m a_l}$ for $a_l\in[0,1]$, one has that
\begin{align}
\Pi_{l=0}^{k-1}c_l&\leq e^{-\mu(2p_{min}-1)\sum_{l=0}^{k-1}\alpha_l}       \nonumber\\
&\leq\frac{1}{(k+1)^3},            \label{new3}
\end{align}
and
\begin{align}
\sum_{m=1}^{k-1}e^{-\sum_{l=m}^{k-1}\frac{3}{l+1}}\frac{1}{m^2}&\leq \sum_{m=1}^{k-1}\frac{(m+1)^3}{(k+1)^3}\frac{1}{m^2}       \nonumber\\
&\leq \sum_{m=1}^{k-1}\frac{8m^3}{(k+1)^3m^2}                   \nonumber\\
&\leq \frac{4}{k}.                               \label{new4}
\end{align}

Then inserting (\ref{new3}) and (\ref{new4}) to (\ref{new2}) leads to the conclusion (\ref{27}). The proof is complete.
\qed

\subsection{Proof of Theorem \ref{t4}}

To ease the exposition, define $\hat{g}_k:=\{\tilde{g}_k^m,m\in[M]\}$. Invoking Assumptions \ref{a3} and Algorithm \ref{al3} yields
\begin{align}
V_{k+1}-V_k\leq-\alpha M_k g_k^\top \operatorname{sign}(\hat{g}_k^s)+\frac{L\alpha^2}{2}M_k^2,                \nonumber
\end{align}
which, by taking the conditional expectation, implies that
\begin{align}
\mathbb{E}(V_{k+1}|x_k)-V_k&\leq -\alpha g_k^\top \mathbb{E}(M_k\operatorname{sign}(\hat{g}_k^s)|x_k)     \nonumber\\
&\hspace{0.4cm}+\frac{L\alpha^2}{2}\mathbb{E}(M_k^2|x_k).                \label{pf17}
\end{align}
For $g_k^\top\mathbb{E}(M_k\operatorname{sign}(\hat{g}_k^s)|x_k)$ in (\ref{pf17}), one has
\begin{align}
g_k^\top\mathbb{E}(M_k\operatorname{sign}(\hat{g}_k^s)|x_k)&=g_k^\top\mathbb{E}(M_k\mathbb{E}(\operatorname{sign}(\hat{g}_k^s)|\hat{g}_k)|x_k)      \nonumber\\
&\hspace{-1.6cm}=\mathbb{E}(M_k\sum_{i=1}^d g_{k,i}\mathbb{E}(\operatorname{sign}(\hat{g}_{k,i}^s)|\hat{g}_k)|x_k)      \nonumber\\
&\hspace{-1.6cm}=\mathbb{E}(M_k\sum_{i=1}^d |g_{k,i}|\mathbb{E}(\operatorname{sign}(\hat{g}_{k,i}^s g_{k,i})|\hat{g}_k)|x_k)      \nonumber\\
&\hspace{-1.6cm}=\mathbb{E}(M_k\sum_{i=1}^d |g_{k,i}|(2I_{p_{k,i}}(\kappa,\kappa)-1)|x_k)         \nonumber\\
&\hspace{-1.6cm}\geq (2I_{p_{min}}(\kappa,\kappa)-1)\|g_k\|_1\mathbb{E}(M_k|x_k),           \label{pf18}
\end{align}
where the last equality has exploited Lemma 7 in \cite{safaryan2019stochastic}, and the inequality comes from the fact that $p_{k,i}\geq p_{min}$ for $i\in[d],k\geq 0$.

As for the last term in (\ref{pf18}), it can be concluded that
\begin{align}
\mathbb{E}(M_k|x_k)&\geq\mathbb{E}(\big\|\frac{1}{M}\sum_{m=1}^M\tilde{g}_k^m\big\|_1|x_k)        \nonumber\\
&\geq \big\|\mathbb{E}(\frac{1}{M}\sum_{m=1}^M\tilde{g}_k^m|x_k)\big\|_1                        \nonumber\\
&=\|g_k\|_1,             \nonumber
\end{align}
which, combining with (\ref{pf17}) and (\ref{pf18}), leads to
\begin{align}
\mathbb{E}(V_{k+1})-\mathbb{E}(V_k)&\leq -\alpha (2I_{p_{min}}(\kappa,\kappa)-1)\mathbb{E}(\|g_k\|_1^2)     \nonumber\\
&\hspace{0.4cm}+\frac{L\alpha^2}{2}\mathbb{E}(M_k^2).                \label{pf19}
\end{align}

Now, for the last term in (\ref{pf19}), one has
\begin{align}
\mathbb{E}(M_k^2)&\leq \frac{1}{M}\sum_{m=1}^M \mathbb{E}(\|\tilde{g}_k^m\|_1^2)          \nonumber\\
&\leq \frac{2}{M}\sum_{m=1}^M \mathbb{E}(\|\tilde{g}_k^m-g_k\|_1^2)+2\mathbb{E}(\|g_k\|_1^2)       \nonumber\\
&\leq 2\sigma^2+2\mathbb{E}(\|g_k\|_1^2),                 \nonumber
\end{align}
which, together with (\ref{pf19}), yields
\begin{align}
&\mathbb{E}(V_{k+1})-\mathbb{E}(V_k)          \nonumber\\
&\hspace{0.4cm}\leq -\alpha (2I_{p_{min}}(\kappa,\kappa)-1-L\alpha)\mathbb{E}(\|g_k\|_1^2)+L\sigma^2\alpha^2.     \nonumber
\end{align}
The rest of proof is similar to that after (\ref{pf16}). This ends the proof.
\qed




\end{document}